%% file: all.tex
\input defv7.tex


\input 1int.tex

\null
\vfill
\eject


\input 2pro.tex


\null
\vfill
\eject

\input 3bib.tex

\end

%% file: defv7.tex
 
\def\a{\alpha}
\def\b{\beta}

\def\g{\gamma}

\def\O{\Omega}
\def\o{\omega}

\def\s{\sigma}

\def \dd#1{{\bf#1}}

\def\cl#1{{\cal#1}}



\def\ouv#1{\smash{\mathop{#1}\limits^{\lower 1pt\hbox
{$\scriptscriptstyle\circ$}}}}

\def\hfl#1#2{\smash{\mathop{\hbox to 12mm{\rightarrowfill}}
\limits^{\scriptstyle#1}_{\scriptstyle#2}}}


\long\def\eno#1#2{\par\smallskip{\bf{#1}}{\it\ {#2}}\par\medskip}

\def\tit#1{\vskip 5mm plus 1mm minus 2mm {\tir #1}
		\vskip 3mm plus 1mm minus 2mm}

\def\stit#1{\vskip 3mm plus 1mm minus 2mm {\bf{#1}}
		\smallskip}

\long\def\dem#1#2{{\bf {#1}}{\ {#2}$\diamondsuit$}\medskip}

\font\tir=cmbx10 at 12pt

\def\ref#1#2#3#4{{\bf #1}{\ #2}{\it ,\ #3}{,\ #4}\medskip}


\def \picture #1 by #2 (#3){\midinsert \centerline 
{\vbox to #2{\hrule width #1 heigth 0pt 
depth 0pt \null \vfill \special {picture #3}}}\endinsert }

\def\scaledpicture #1 by #2 (#3 scaled #4) {{
\dimen0 =#1 \dimen1 =$2
\divide \dimen0 by 1000 \multiply \dimen0 by #4
\divide \dimen1 by 1000 \multiply \dimen1 by #4
\picture \dimen0 by \dimen1 (#3 scaled $4)}}

\def\figure #1 #2 #3 {\midinsert \vglue 3mm 
{\vbox to #3 {\hrule width 6cm height 0cm depth 0cm \vfill
{\special {picture #1 scaled #2}}}}\vglue 2mm \endinsert}

\magnification=1200

%% file: 1int.tex
\input defv7.tex
\overfullrule=0pt

{\centerline {\tir {Convergence or generic divergence of }}}

{\centerline {\tir {Birkhoff normal form}}}

\bigskip
\bigskip

{\centerline {{\bf Ricardo Perez-Marco}\footnote {*} 
{  UCLA, 
Dept. of Mathematics, 
405, Hilgard Ave., Los Angeles, CA-90095-1555, 
USA, e-mail: ricardo@math.ucla.edu. 
}}}

\bigskip

{\bf Abstract.} {\it 
We prove that Birkhoff normal form of hamiltonian flows
at a non-resonant singular point with given quadratic part are 
always convergent or generically divergent. The same result is 
proved for the normalization mapping and any formal first integral.
}

\bigskip


\medskip


\bigskip
\bigskip

\bigskip

{\centerline  {\bf Introduction }}

\bigskip
\bigskip

In this article we study analytic ($\dd R$ or $\dd C$-analytic) 
hamiltonian flows
$$
\eqalign {
\dot x_k &=+{\partial H \over \partial y_k} \cr
\dot y_k &=-{\partial H \over \partial x_k}\cr
}
$$
where $x_k, y_k \in \dd C$ (resp. $\dd R$), $k=1,2, \ldots n$, and 
$H$ is an analytic hamiltonian with power
series expansion at $0$ beginning with quadratic terms
(so $0$ is a singular point of the analytic vector
field).
We shall restrict our attention to those $H$ having
a non-resonant quadratic parts: If $(\lambda_1 , \ldots ,
\lambda_n )$ are the 
eigenvalues of the symmetric matrix $Q$ where 
${1\over 2} (x,y)Q (x,y)^t$ is the quadratic part of $H$
then, defining $\lambda_{n+1}=-\lambda_1 , \ldots 
\lambda_{2n}=-\lambda_n$, there is no relation of the form
$$
i_1\lambda_1 +\ldots +i_{2n} \lambda_{2n} =0
$$ 
with integral coefficients $i_1 , \ldots , i_{2n}$
except for the trivial
case $i_1=\ldots =i_{2n} =0$. Due to some confusion that 
one finds in some of the litterature on the problem 
of convergence of Birkhoff normal 
form and Birkhoff transformation, we start with 
a brief historical overview.

\bigskip

The normal form of a hamiltonian flow near 
a singular point has been studied since the 
origins of mechanics. The long time evolution of the 
system near the equilibrium position is better controlled
in variables oscullating those of the normal form that 
corresponds to a completely integrable system.
This idea is at the base of many computations in 
Celestial Mechanics. Its importance, both practical 
and theoretical, cannot be underestimated.
One can consult the reference
memoir "Les m\'ethodes nouvelles de la m\'ecanique 
c\'eleste" by H. Poincar\'e ([Po]) to get an idea of the 
central place that the perturbative approach played already 
in the XIXth century.

\medskip

Assuming that the eigenvalues of the quadratic part
of $H$ present no resonances, we have a simple formal 
normal form. This result goes back to C.E. Delaunay [De] and 
A. Lindstedt [Li] (see [Po], [Si2]).
Nowadays this normal form is named after Birkhoff. 
Birkhoff normal form is the starting point of most of 
the studies of stability near the equilibrium point: 
the first studies by E.T. Whittaker
[Wh],
T.M. Cherry [Ch], G.D. Birkhoff [Bi1, [Bi2], 
and C.L. Siegel [Si1] [Si2], 
K.A.M. theory ([Ko], [Ar], [Mo]), 
Nehoroshev's diffusion estimates [Ne],...

\medskip

The dream of an analytic conjugacy to the normal form 
(without fixing the quadratic part of $H$) was quickly 
dissipated after the work of H. Poincar\'e ([Po] vol.I chapitre V). 
Poincar\'e's divergence theorem is the starting point of 
his difficult proof of the inexistence of non-trivial 
local first integrals
in the three body problem.

Research then 
focused in understanding the divergence of the conjugation mapping
(normalization mapping) with a fixed non-resonant quadratic part for $H$. The normal form is unique. 
The normalization mapping is not unique, but 
appropriate normalizations determine it uniquely.
Different results showed with increasing
strength that the normalization mapping was generically divergent.
We refer to the book of C.L. Siegel and J. Moser  
([Si-Mo] chapter 30) for an overview. The strongest result
on divergence being proved by Siegel in 1954 ([Si2]) showing 
the generic divergence of the normalization, the quadratic part 
of the hamiltonian being fixed but otherwise arbitrary. 
A.D. Bruno [Br] and 
H. R\"ussman [Ru2] [Ru3] proved the convergence of the normalization
when Birkhoff normal form for the hamiltonian is quadratic 
and the eigenvalues satisfy
Bruno's arithmetic condition (other proofs can be found in [El2], 
[E-V], [Sto1], [Sto2]).

\medskip

Despite this progress, the most natural question
remains untouched. The question  
is not the convergence or divergence
of the normalizing map, but actually the convergence or divergence
of Birkhoff normal form itself. 
If in first place Birkhoff
normal form is diverging, then there is 
no point in trying to conjugate
to the normal form. Also in this case the normalization 
is necessarily diverging.

\medskip

Very surprisingly, there seems to be no significant results on this
fundamental question. It appears to be a very hard question.
The author first learned about it
from H. Eliasson.
The references in the literature are 
scarce. H. Eliasson points out in the introduction of his article [El1]
that  

\medskip

{\it "...if the normal form itself  is convergent or divergent 
is not known..."}, 

\medskip

and he points out in [El2] 

\medskip

{\it "...Generically (...) 
the formal transformation is divergent.
(if the normal form itself also is generically divergent is not known)."}

\medskip

These are the only citations in the literature that the 
author is aware of (despite the title of [It] what is really
proved there is the convergence of the normalization).
On the other hand, one frequently finds in some litterature
the claim "Birkhoff normal form is generically diverging" in place 
of the "Birkhoff transformation is generically diverging"...

\bigskip

More surprisingly, not a single example is known of an analytic hamiltonian
having a divergent Birkhoff normal form.
The main result in this article is that 
the existence of a single example 
with divergent Birkhoff normal form forces 
generic divergence. 
To be more precise we need to introduce the notion of $\Gamma$-capacity
of a subset of $\dd C^n$. This notion generalizes the notion 
of capacity in dimension $1$. We recall the definition in section
2. We refer the reader to [Ro].
An important property, as in dimension $1$,  
is that a set $E\subset \dd C^n$ with zero 
$\Gamma$-capacity is Lebesgue and Baire thin, i.e. $E$ 
has zero Lebesgue measure and is of the first category (a countable 
union of nowhere dense sets).

In order to talk about generic properties we define a 
natural Baire 
space. We consider the Banach space $\cl H$ of Hamiltonians 
with radius of convergence $1$ endowed with a uniform norm
in some open subset of the disk of convergence.
Similar results hold for $\dd C$-analytic
and $\dd R$-analytic hamiltonians.

We can now state:

\eno {Theorem 1.}{ We consider the subspace of $\cl H_Q\subset \cl H $ 
of analytic  hamiltonians 
$$
H=\sum_{l=2}^{+\infty } H_l
$$ 
with fixed non-resonant quadratic part $H_2$ given by the symmetric 
matrix $Q$.

If there exists one hamiltonian $H_0 \in \cl H_Q$ with 
divergent Birkhoff normal form (resp. normalization), then a generic 
hamiltonian in $\cl H_Q$ has divergent 
Birkhoff normal form (resp. normalization). 

More precisely, all hamiltonians in 
any complex (resp. real) affine finite dimensional subspace $V$ of 
$\cl H_Q$ have a convergent Birkhoff normal form (or normalization), 
or only an
exceptional subset in $V$ of $\Gamma$-capacity $0$ (resp. of Lebesgue 
measure $0$) have this property.
}

Observe that the second scenario holds for all affine subspaces
containing $H_0$. The result obtained in the real analytic case
is stronger than stated. When $V$ is a one real dimensional affine 
line, the exceptional set has zero capacity in the complexification
of $V$. So the exceptional set has even Hausdorff dimension zero.

The important issue that remains unsettled is thus the existence of 
hamiltonians with diverging Birkhoff normal form for any non-resonant
quadratic part. It seems 
to be the prevalent opinion among specialist that there is 
generic divergence for all non-resonant quadratic 
parts. This feeling is probably motivated by the 
divergence results on the normalization, which, it is worth 
noting, are independent of the quadratic part. The author 
sees no reason against
the convergence of Birkhoff normal forms, in particular when
the eigenvalues of the quadratic part of $H$ enjoy good arithmetic
properties. Fixing the quadratic part of the hamiltonian, the answer
may depend on the arithmetic of its eigenvalues.

\bigskip

On the other hand, using standard methods of Small Divisors, it is 
not difficult to exhibit hamiltonians with diverging normalizations
using Liouville eigenvalues for the quadratic part.
Combining this construction with the previous theorem, one recovers
with a simple proof Siegel's result ([Si2]) on the generic divergence of the 
normalization mapping for some fixed quadratic parts.

\bigskip

Note that fixing the quadratic part of the hamiltonian makes
the problem much harder, not allowing to take any advantage of the 
arithmetic of the eigenvalues. One can find
in the litterature results whithout fixing the quadratic part ([Po]
volume I chapter V, [Koz]). One may ask about the reason for  
studying hamiltonians with fixed quadratic part. Note that for 
systems with particles, the masses enter directly into the quadratic
part of the hamiltonian through the kinetic 
energy. Thus if one, for example, wants to show the non-integrability
of a given system {\it with given masses} then families of hamiltonians
with fixed quadratic part arise naturally. One can cite at this respect
the strict criticism 
of A. Wintner to Poincar\'e's proof of non-integrability
of the three body problem ([Wi] 320, p.241):

\medskip

{\it ...Poincar\'e has stablished a result which concerns the 
non-existence of additional integrals (...) Nevertheless, his result,
as well as its formal refinement obtained by Painlev\'e, is not 
satisfactory (...) In fact, these negative results do not deal 
with the case of fixed, but rather with unspecified, values of the 
masses $m_i$ (...) Clearly, these assumptions in themselves do not 
allow any dynamical interpretation, since a dynamical system is 
determined by a fixed set of positive numbers $m_i$ ...}

\medskip

Whitout sharing this extreme view, one cannot deny some 
point in Wintner's criticism.

\bigskip

We prove a second theorem on the divergence of first integrals.
The classical approach to integrability of hamiltonian systems
is based on {\it first integrals}. A first integral $P$ is a 
convergent power series in the $2n$ variables $x_1, \ldots , y_n$
such that 
$$
\{ P, H \}=0
$$
where the Poisson bracket is defined by 
$$
\{ P , H \} =\sum_{k=1}^n \left ( {\partial P \over \partial 
x_k }{\partial H \over \partial y_k} -{\partial P \over \partial 
y_k}{\partial H \over \partial x_k} \right ) \ .
$$
The equation $\{P, H \}=0$ is equivalent to $\dot P =0$, that is
to the conservation of $P$. By E. Noether's theorem, symmetries 
of the hamiltonian generate first integrals. Two first integrals, 
$P_1$ and $P_2$, are in involution (or functionally independent) 
if their Poisson 
bracket vanishes
$$
\{ P_1 , P_2 \}=0 \ .
$$ 
At a non-singular point of the hamiltonian, Liouville's theorem shows
that the hamiltonian system is integrable
by quadratures if there exists $n$ first integrals in 
involution. 
The case of a non-resonant singular point as considered here is 
more involved. It has been shown by H. R\"ussman [Ru1] for 
$n=2$ and in general by J. Vey [Ve] and H. Ito [It]
that the existence of $n$
first integrals in involution forces the convergence of the 
normalization to Birkhoff normal form (H. Eliasson settled 
the analogue of Vey's theorem in the 
$C^\infty$ case [El1], [El3]). Recently L. Stolovitch 
found a unified approach 
to Bruno's theorem cited before and Vey's and Ito's theorems 
([Sto1], [Sto2]).
Once all symmetries of a system have been used to find 
first integrals in involution, the natural question is if 
there are any others. 
Multiple approaches to non integrability 
have been developped starting from H. Poincar\'e.
We refer to 
[Koz] for an overview of classical methods. 
R. de la Llave has recently found that Poincar\'e's
conditions are necessary and sufficient for uniform 
integrability ([Ll], see also the paper by G. 
Gallavotti [Ga]). 
We refer to [Mo] for an account on 
recent methods of S.L. Ziglin, J. Morales Ruiz and 
J.-P. Ramis. In the smooth non-analytic setting we refer
to the work of R.C. Robinson ([Rob]).

It is natural to define 
the {\it degree of integrability} of a hamiltonian as the maximal 
number $1\leq \iota (H) \leq n$ of first integrals in involution.
When the normalization is convergent, we have 
that $\iota (H)=n$, so the study of convergent first integrals can 
be seen as a refinement of the study of the convergence of 
the normalization.

\eno {Theorem 2.}{We consider the space $\cl H_Q$.
Given a hamiltonian $H_0\in \cl H_Q$, we have for a generic 
hamiltonian $H\in \cl H_Q$,
$$
\iota (H) \leq \iota (H_0) \ .
$$
More precisely, let $P$ be a universal formal first integral. In any 
complex (resp. real) affine finite dimensional subspace $V$
of $\cl H_Q$ all hamiltonians $H\in V$ have converging $P(H)$,
or only an exceptional set in $V$ of $\Gamma$-capacity zero (resp.
Lebesgue measure zero) have this property.}

We give in section 1 a precise definition of universal 
formal first integral.
This theorem reduces the proof of the generic divergence of a given 
formal first integral in a family of hamiltonians, 
to the divergence for one hamiltonian.
Also, given a family $V$, the minimum degree of integrability in $V$,
$$
\iota_V =\min_{H\in V} \iota (H)
$$
is attained for a generic $H\in V$.

\bigskip

The families $V$ in theorem 1 and 2 can be more general than finite
dimensional affine subspaces. The same proof gives the results
for example when $V$ is parametrized polynomially by $\dd C^m$.
It is interesting to note how in these theorems the complexification 
of the problem sheds new light on the real analytic case.

The main idea of this article has also been applied 
to other problems of small divisors ([PM1], [PM2]).

\bigskip

{\bf Acknowledgements.}

\medskip

The author is grateful to A. Chenciner and H. Eliasson 
for conversations on the subject.

%% file: 2pro.tex
\input defv7.tex
\tit {1) Birkhoff normal form and first integrals.}

\stit {a) Birkhoff normal form.}

We review briefly in this section the construction of Birkhoff 
normal form.  We follow [Si-Mo]. We need to pay particular 
attention on the polynomial dependence of the transformation 
and Birkhoff normal form on the original coefficients of 
the hamiltonian function. More precisely, it is important 
for our purposes to keep 
track of the degrees of the polynomial dependence.
We use the sub-index notation for partial derivatives.

We consider an analytic hamiltonian ($\dd R$ or $\dd C$ analytic)
$$
H (x,y) =\sum_{l=2}^{+\infty } H_l (x,y)
$$
where $H_l$ is the homogeneous part of degree $l$ in the 
real or complex variables $x_1 ,  \ldots , x_n $,
$ y_1 , \ldots , y_n$.
We can assume, by means of a preliminary linear change of 
variables, that $H_2$ is already in diagonal form ([Bi]
section III.7)
$$
H_2 (x,y)=\sum_{k=1}^n \lambda_k x_k y_k \ .
$$
We look for a simpler normal form of the system
$$
\eqalign {
\dot x_k =&H_{y_k} \cr
\dot y_k =& -H_{x_k}\cr
}
$$
We consider symplectic transformations that leave unchanged
the hamiltonian character of the system of differential 
equations.
The new variables $(\xi , \eta )$ are related to the old 
ones $(x , y)$ by the canonical transformation
$$
\eqalign {
x_k &= \varphi_k (\xi , \eta )=\xi_k +\sum_{l=2}^{+\infty }
\varphi_{kl} (\xi , \eta )\cr
y_k &= \psi_k (\xi , \eta )=\eta_k +\sum_{l=2}^{+\infty } 
\psi_{kl} (\xi , \eta )\cr
}
$$
where the $\varphi_{kl}$ and $\psi_{kl}$ are the homogeneous
parts of degree $l$.
These canonical transformations are defined by a generating 
function 
$$
v(x,\eta)=\sum_{l=2}^{+\infty } v_l (x, \eta )
$$
where $v_l$ is the homogeneous part of degree $l$, 
and $v_2 (x , \eta )=\sum_{k=1}^{+\infty } x_k \eta_k $. Then 
the canonical transformation is defined by the equations
$$
\eqalign {
\xi_k &=v_{\eta_k }(x, \eta )=x_k 
+\sum_{l=3}^{+\infty } v_{l, \eta_k} (x, \eta ) \cr
y_k &=v_{x_k} (x , \eta )=\eta_k +\sum_{l=3}^{+\infty } v_{l, x_k } (x, \eta )\cr
}
$$
So we get 
$$
\eqalign {
x_k &=\xi_k 
-\sum_{l=3}^{+\infty } v_{l \eta_k} (\varphi (\xi , \eta), \eta ) \cr
y_k &=\eta_k +\sum_{l=3}^{+\infty } v_{l x_k } (\varphi (\xi , \eta ), \eta )\cr
}
$$
and
$$
\eqalign {
\varphi_{kl} (\xi , \eta ) &=-v_{l+1 ,\eta_k} (\xi , \eta )
-\left \{ \sum_{j=3}^l v_{j, \eta_k} (\varphi (\xi ,\eta ) , \eta ) \right \}_l
\cr
\psi_{kl} (\xi , \eta ) &=v_{l+1 ,x_k} (\xi , \eta )
+\left \{ \sum_{j=3}^l v_{j, x_k} (\varphi (\xi ,\eta ) , \eta ) \right \}_l
\cr
}
$$
where $\{ . \}_l$ indicates that we take the $l$ homogeneous part
of the expression within brackets. From these expressions we have that 
the coefficients of $\varphi_{kl}$ and $\psi_{kl}$ are polynomials 
with integer coefficients on the coefficients of $v_3 , \ldots , v_l ,
v_{l+1}$.

To each coefficient of $v_l$ we assign a degree $l-2$ (as we will 
see next, we will choose a canonical transformation so that 
the coefficients of the $v_l$'s are polynomials on the coefficients
of $H$ of degree $l-2$ at most). By induction, we show that 
the degree of $\varphi_{kl}$ is at most $l-1$. For $l=2$ it is clear.
Then by induction, the degree of the coefficients of  
the homogeneous part of degree $l$ of an homogeneous monomial 
$$
\prod_{k=1}^n (\varphi (\xi , \eta ))^{\a_k } \eta^{\b_k}
$$
of total degree $j$ ($\sum \a_k +\sum \b_k =j$)
is at most $l-j$. Thus the degree of 
$$
v_{j, \eta_k} (\varphi (\xi ,\eta ) , \eta )
$$
is at most $(j-2)+(l-j+1)=l-1$, and this finishes the induction.
The same discussion applies to $\psi$ and the coefficient 
$\psi_{kl}$ has degree $l-1$.

\bigskip

Now the canonical transformation generated by $v$ transforms 
the differential system into 
$$
\eqalign {
\dot \xi_k &=K_{\eta_k } \cr
\dot \eta_k &= -K_{\xi_k} \cr
}
$$
where 
$$
K(\xi , \eta )=\sum_{l=2}^{+\infty } H_l (\varphi (\xi , \eta) , 
\psi (\xi , \eta ) ) =\sum_{l=2}^{+\infty } K_l (\xi , \eta )
$$
where  $K_l$ is the $l$-homogeneous part.

Our aim is to construct a canonical transformation which gives 
a hamiltonian $K$ only depending on power series of the 
products $\omega_k =\xi_k \eta_k$.
The coefficients of $v$ are constructed by induction on the degree
$l$ of the homogeneous part. Assume that the choices for $v_3 , 
\ldots , v_{l-1}$  have been done so that the new hamiltonian
has monomials of order $\leq l-1$ only depending on the $\omega_k$'s.
We consider a monomial of degree $l$
$$
P=\prod_{k=1}^n \xi_k^{\a_k } \eta_k^{\b_k } \ .
$$
We want to choose the coefficient $\g$ of $P$ in $v_l(\varphi (\xi , \eta),
\eta )$ such 
that the new hamiltonian does not contain the monomial $P$.
Note that 
$$
K_l (\xi , \eta )=\sum_{k=1}^{+\infty } \lambda_k \left (
\xi_k v_{lx_k}(\varphi (\xi ,\eta ), \eta) -\eta_k v_{l \eta_k} (
\varphi (\xi ,\eta ), \eta) \right ) + A
$$
where the first term comes from the expansion of 
$H_2 (\phi (\xi , \eta ), \psi (\xi , \eta))$ and the second term $A$
collects everything coming from higher order. The coefficients 
in the expression $A$ 
are polynomials in the coefficients of $v_3 , \ldots , v_{l-1}$
and linear functions in the coefficients of $H_3 , \ldots H_l$.

By induction we prove at the same time that the coefficients of
$v_l$ are polynomials of degree $l-2$ on the coefficients of 
$H_3 , \ldots , H_l$, and also the coefficients of $K_l$ are 
polynomials of degree $l-2$ on the coefficients of 
$H_3 ,\ldots , H_l$. Assuming the induction hypothesis, 
we have as before that the right hand side  in the above formula
for $K_l$ is a polynomial of degree $\leq l-2$ on the coefficients
of $H_3 , \ldots , H_l$.

Now we have 
$$
\sum_{k=1}^n \lambda_k (\xi_k P_{\xi_k} -\eta_k P_{\eta_k})=
 \left ( \sum_{k=1}^n \lambda_k (\a_k -\b_k ) \right ) P
$$
Thus if $\lambda = \sum_{k=1}^n \lambda_k (\a_k -\b_k )\not= 0$,
choosing 
$$
\g =- {1\over \lambda} \left \{ A \right \}_{P}
$$
(where brackets indicate that we extract the $P$ monomial)
the new hamiltonian will not contain the monomial $P$.
Note that by the non-resonance condition, $\lambda =0$ only 
happens when $\a_k =\b_k$ for $k=1, \ldots , n$.  In that way
we determine all coefficients of $v_l$ except those of the 
monomials which are a product of $\omega_k$'s. Note also that
by induction these coefficients are polynomials on the 
coefficients of $H_3, \ldots , H_l$ of degree $\leq l-2$.

In order to determine the coefficients of $v_l$ for the 
remaining monomials one takes the normalization that no 
product of powers of $\omega_k$'s appears in 
$$
\Phi =\sum_{k=1}^n (\xi_k y_k -\eta_k x_k) 
$$
when expressed in $(\xi , \eta )$ variables.
One checks that this determines uniquely $v$ and thus the 
canonical transformation that transforms the hamiltonian
into its Birkhoff normal form. 
When $H$ is real analytic, it is easy to check ([Si-Mo]) that
the previous construction yields a real formal canonical 
transformation and a real Birkhoff normal form.
We summarize this discussion
in the following proposition.

\eno {Proposition 1.1.}{Given a hamiltonion flow 
$$
\eqalign {
\dot x_k &=H_{y_k} \cr
\dot y_k &=-H_{x_k} \cr
}
$$
with $H(x,y)=\sum_{l=2}^{+\infty } H_l (x,y)$ with non-resonant
quadratic part $H_2$,
there exists a unique formal canonical transformation defined
by a formal generating series
$$
v(x ,\eta )=\sum_{l=2}^{+\infty } v_l (x, \eta )
$$
such that in the new variables $(\xi_k , \eta_k)$ the differential
system takes the form
$$
\eqalign {
\dot \xi_k &=K_{\eta_k} \cr
\dot \eta_k &=-K_{\xi_k} \cr
}
$$
where the new hamiltonian $K$ is a formal power series in the 
products $\omega_k =\xi_k \eta_k$, and the expression
$$
\Phi =\sum_{k=1}^n (\xi_k y_k -\eta_k x_k) \ .
$$
contains no product of the $\omega_k$ in the $(\xi , \eta )$ variables.
Moreover, the coefficients of the homogeneous part of $K$ of 
degree $l$ and of $v_l$ 
are polynomials of degree $l-2$ in the coefficients 
of $H_3 ,\ldots , H_l$.
}

\stit {b) First integrals.}

We review some classical facts about first integrals (see [Si1]).

If the normalization is converging, then all expresions 
$$
\omega_k =\xi_k \eta_k
$$
are first integrals since
$$
\{ \omega_k , K \} =\eta_k K_{\eta_k}-\eta_k K_{\xi_k}
=\xi_k \eta_k (K'-K') =0 \ .
$$
Expressing $\omega_k$ in terms of the initial variables $(x,y)$
we get $n$ formal first integrals
$$
P_k (x, y)=\xi_k (x,y) \eta_k (x,y)\ .
$$
Observe that 
$$
\eta_k =y_k -\sum_{l=3}^{+\infty } v_{l, x_k} (x, \eta) \ .
$$
So if 
$$
\eta_k (x,y)=y_k +\sum_{l=2}^{+\infty } \eta_{kl} (x,y)
$$
where $\eta_{kl}$ is the $l$-homogeneous part of $\eta$, then
by induction the coefficients of $\eta_{kl}$ are polynomial
on the coefficients of $H_3 , \ldots , H_{l+1}$ of degree $l-1$.

We reach the same conclusion for $\xi_k$ using
$$
\xi_k (x,y)=v_{\eta_k} (x, \eta) =x_k +\sum_{l=3}^{+\infty } 
v_{l,\eta_k } (x , \eta)  \ .
$$

Now, we have the following formal lemma ([Si1] lemma 1):

\eno {Lemma 1.2.}{Any formal integral $P$ can be represented 
as a formal power series in the $n$ first integrals $\o_1 , 
\ldots , \o_n$.}

\dem {Proof.}{Let $P(x , y)$ be a formal first integral.
We have that
$$
P(x,y)=\hat P (\xi , \eta )
$$
is a formal first integral in the $(\xi , \eta)$ variables.
We write
$$
\hat P=T+J
$$
where $T$ is the formal power series containing all monomials
of the form
$$
\prod_{k=1}^n \xi_k^{\a_k}\eta_k^{\a_k} \ ,
$$
thus $T$ is a formal power series on the $n$ formal first 
integrals $\o_1 , \ldots , \o_n$. We only need to show that 
$F$ is identically $0$.  If not consider the leading monomial
of $J$
$$
L=\prod_{k=1}^n \xi_k^{\a_k}\eta_k^{\b_k} \ ,
$$
with some $\a_k-\b_k \not=0$.
The formal power series $J$ is a formal first integral, and
computing the leading term in
$$
\eqalign {
0=\{ J , K \} &= \sum_{k=1}^{+\infty } J_{\xi_k } K_{\eta_k }-
J_{\eta_k} K_{\xi_k } \cr
&= \sum_{k=1}^{+\infty }(J_{\xi_k }\xi_k -J_{\eta_k } \eta_k )K' \cr
&=\sum_{k=1}^{+\infty } ((\a_k-\b_k ) L+\ldots ) K'
}
$$
we get
$$
\sum_{k=1}^{+\infty } \lambda_k (\a_k -\b_k) =0 \ .
$$
So by the non-resonance condition, 
for $k=1, \ldots , n$, $\a_k -\b_k =0$ and $J=0$.}

Thus we can identify the set of formal first integrals with 
the formal power series in $n$ variables. 

\eno {Definition 1.3.}{A universal formal first integral
$P(H)$ is $P(H)=F(\o_1 , \o_n)$ where $F$ is a formal 
power series in $n$ variables.}

\eno {Corollary 1.4.}{ Any universal 
formal first integral $P(H)$ has
coefficients that are monomials of degree $l$ depending polynomially
on the coefficients of $H_3 , \ldots , H_{l+1}$ with degree 
$\leq l-1$.}

\null
\vfill
\eject

\tit {2) Proof of the theorems.}

\stit {a) Potential theory.}

\stit {$\Gamma$-capacity.}

We recall the definition of $\Gamma$-capacity and 
we refer to [Ro] for more properties.
Let $E \subset \dd C^m$. The $\Gamma$-projection of $E$ on 
$\dd C^{m-1}$ is the set $\Gamma_m^{m-1} (E)$ 
of $z=(z_1, \ldots , z_{m-1})\in \dd C^{m-1}$
such that 
$$
E\cap \{ (z, w)\in \dd C^m \}
$$
has positive capacity in the complex plane 
$ \dd C_z=\{ (z, w)\in \dd C^m \}$.
We define 
$$
\Gamma_m^1 (E)=\Gamma_2^1 \circ \Gamma_3^2\circ \ldots \Gamma_m^{m-1} (E) \ .
$$
Finally, the $\Gamma$-capacity is defined as 
$$
\Gamma {\hbox {\rm -Cap}} (E)=\sup_{A\in U(m,\dd C)} {\hbox {\rm Cap}}
\ \Gamma_m^1 (A(E)) \ .
$$
where $A$ runs over all unitary transformations of $\dd C^m$.

The following lemma is useful ([Ro] Lemma 2.2.8 p.92)

\eno {Lemma.}{Let $E \subset \dd C^m$, $E\not=\dd C^m$ and assume 
that the intersection of $E$ with any complex line which is not 
a subset of $E$ has inner capacity zero. Then the $\Gamma$-capacity 
of $E$ is zero.}

As we will see, the set of elements in $\cl H$ 
with convergent Birkhoff normal 
form (or normalization) is an $F_{\sigma}$-set, so capacitable,
and the inner capacity is the capacity of the set. 
Thus using this lemma, we are reduced to 
prove the second assertion of the theorem 
only when the sub-space $V$ of $\cl H$ ahs dimension $1$.

\stit {Bernstein lemma.}

The following is a classical lemma in potential theory and 
approximation theory ([Ra] p.156). It plays a crucial role in the 
proof of theorem 1.

\eno {Lemma (Bernstein).}{Let $E\subset \dd C$ be a non-polar
compact set (i.e. ${\hbox {\rm cap}} (E) >0$). Let $\O$ be the connected 
component of ${\overline {\dd C}}-E$ containing $\infty$.
Then for any polynomial $P$ of degree $n$, we have for $t\in \dd C$,
$$
|P(t)|\leq e^{n g_{\O} (t,\infty)} \ ||P||_{C^0(K)} 
$$
where $g_{\O}$ denotes the Green function of $\O$.
}

The proof is quite simple, we include it here for completeness.

\dem {Proof.}{We can assume the polynomial monic. Then
$$
u(t)=\log P(t)-\log ||P||_{C^0(K)} -g_{\O } (t, \infty)
$$
is sub-harmonic, is negative near $\infty$ (because 
$g_{\O }(t, \infty )=\log |t| +{\hbox {\rm cap}} (E) +o(1)$),
and $\limsup u(t) \leq 0$ when $t\to K$. The 
application of the maximum principle
concludes the proof.}

\stit {b) Proof of theorem 1.}

The assertion about the divergence of the normalization 
mapping follows the same lines than the case of the 
Birkhoff normal. The convergence or divergence 
of the normalizing transformation is equivalent 
to the convergence or divergence of the generating 
function. Then the proof proceeds in the same way as 
below using the the polynomial dependence 
of the generating function on the coefficients of $H$
(proposition 1.1).

For the elementary construction of hamiltonians with 
divergent normalization mentioned at the end of the introduction,
we refer the reader to the end of section 30 of [Si-Mo], and 
to Siegel's article [Si1].

\medskip

We consider the problem of convergence or divergence 
of Birkhoff normal form.
The first assertion of the theorem follows from the second. 
Actually,
consider the set $F_n \subset \cl H_Q$ of hamiltonians
having a converging Birkhoff normal form
with radius of convergence $> 1/n$, and bounded by $1$
in the ball of radius $1/n$. This set $F_n$ is closed, 
and 
$$
F=\bigcup_{n\geq 1} F_n
$$
is the set of all hamiltonians in $\cl H_Q$ having a 
convergent Birkhoff normal form (so this set is an
$F_\s$-set).
Moreover, the open set
$\cl H_Q-F_n$ is dense. Otherwise  let $H_1$ be a hamiltonian
in the interior of $F_n$. Considering the complex (resp. real) 
affine subspace
$V=\{ (1-t)H_0+tH_1 ; t\in \dd C ({\hbox {\rm resp. }} \dd R )\} 
\subset H_Q$ we have, 
according to the second assertion in theorem 1, 
that the set of hamiltonians
with converging Birkhoff normal form  must have capacity zero
(resp. Lebesgue measure $0$).
But on the other hand it contains a neighborhood of $1$. 
Contradiction.

\bigskip
The real analytic result follows from the $\dd C$-analytic 
one by the observation that the intersection of a set of 
$\Gamma$-capacity $0$ in $\dd C^n$ with $\dd R^n \subset \dd C^n$
has Lebesgue measure $0$ (see [Ro] Lemma 2.2.7 p. 90).
 
We consider a complex finite dimensional affine subspace $V$ of $\cl H $.
According to the definition of $\Gamma$-capacity we are reduced
to the case of a one dimensional subspace $V\approx \dd C$.
We can parametrize linearly the coefficients of hamiltonians $H\in V$ with 
a complex parameter $t\in \dd C$, and we denote $H_t$ the 
corresponding hamiltonian in $V$. Note that the coefficients
of $H_t$ are linear functions of $t$.

 We assume that the Birkhoff normal form of hamiltonians $H_t$ 
corresponding to a set of values $t\in E\subset \dd C$ 
of positive capacity (non-polar)
are converging. We want to prove that all the other hamiltonians in $V$ 
have converging Birkhoff normal form. 

We have 
$$
F=\bigcup_{n\geq 1} F_n
$$
where $F_n$ the set of parameters $t\in \dd C$ such that 
the hamiltonian $H_t$ has a Birkhoff normal form  $K_t$ with radius
of convergence larger or equal to $1/n$ and $K_t$ is bounded 
by $1$ in this ball. So if $F$ is non-polar, 
we have  for some 
$n\geq 1$ that $F_n$ is not polar (and this set is also 
closed). If we denote
$$
K_t (\xi , \eta )=\sum_i K_i(t) (\xi , \eta)^i \ ,
$$
then, according to proposition 1.1,  
the coefficients $K_i(t)$ depend polynomially on $t$ with 
degree $\leq |i| -2$ (for $|i|\geq 3$).
Now, there exists $\rho_0 >0$
such that for all $t\in F_n$, 
$$
\varphi (t)=\limsup_{|i|\to +\infty} |K_{i}(t)|\rho_0^{-|i|} <+\infty \ .
$$
The function $\varphi$ is lower semicontinuous, and
$$
F_n=\bigcup_m L_m
$$
where $L_m=\{ z\in F_n ; \varphi (t) \leq m \}$ is closed.  
By Baire theorem for some $p$, $L_m$ has non-empty interior (with 
respect to $F_n$), thus this $L_m$ has positive capacity. Finally
we found a compact set $C=L_m$ of positive capacity such that 
there exists $\rho_1 >0$ such that for any $t\in C$ and 
and all $i\in \dd N^n$,
$$
|K_{i} (t)| \leq \rho_1^{|i|} \ .
$$
Using Bernstein's lemma and proposition 1.1 we get that for any compact 
set $C_0\subset \dd C$ we have for $|i|\geq 3$,
$$
||K_i||_{C^0 (C_0)} \leq \rho^{ |i|-2} \rho_1^{|i|} \ ,
$$
for some constant $\rho$ depending only on $C_0$.
Thus $K_t$ is converging for any $t\in \dd C$.

\stit {c) Proof of theorem 2.}

The proof of theorem 2 goes along the same lines than the 
proof of theorem 1, using the polynomial dependence of  
universal formal first integrals proved in corollary 1.4.

%% file: 3bib.tex
\input defv7.tex
\bigskip

{\centerline {\bf BIBLIOGRAPHY}}

\bigskip
\bigskip

\ref{[Ar]}{V.I. ARNOLD}{Proof of A.N. Kolmogorov's theorem
on the preservation of quasi-periodic motions under small perturbations
of the hamiltonian}
{Usp. Mat. Nauk., {\bf 18}, {\bf 5}, 113, 1963, p.13-40.}

\ref{[Bi1]}{G.D. BIRKHOFF}{Dynamical systems}
{American Mathematical Society, 1927.}

\ref{[Bi2]}{G.D. BIRKHOFF}{Surface transformations and their 
dynamical applications}
{Acta Math. {\bf 43}, 1922, p.1-119.}

\ref{[Br]}{A.D. BRJUNO}{Analytical forms of differential equations}
{Trans. Mosc. Math. Soc., {\bf 25}, 1971, {\bf 26}, 1972.}

\ref{[Ch]}{T.M. CHERRY}{On the solution of hamiltonian systems of 
differential equations in the neighborhood of a singular point}
{Proc. London Math. Soc., II {\bf 27}, 1928, p. 151-170.}

\ref{[De]}{C.E. DELAUNAY}{Th\'eorie du mouvement de la lune}
{Paris Mem. pr\'es., {\bf 28}, 1860, {\bf 29}, 1867.}

\ref{[El1]}{L.H. ELIASSON}{Normal forms for Hamiltonian systems with 
Poisson commuting integrals-elliptic case}
{Comment. Math. Helvetici, {\bf 65}, 1990, p. 4-35.}

\ref{[El2]}{L.H. ELIASSON}{Hamiltonian systems with linear 
form near an invariant torus}
{Non-linear Dynamics, Bologna, 1988, World Sci. Publishing,
Teaneck, NJ, 1989, p.11-29.}

\ref{[El3]}{L.H. ELIASSON}{Hamiltonian systems 
with Poisson commuting integrals}
{Thesis, Univ. Stockholm, 1984.}

\ref{[E-V]}{J. ECALLE, B. VALLET}{Correction and linearization 
of resonant vector fields and diffeomorphisms}
{Math. Z., {\bf 229}, 1998, p.249-318.}

\ref{[Ga]}{G. GALLAVOTTI}{A criterion of integrability of 
perturbed nonresonant harmonic oscillators. "Wick ordering"
of the perturbations in classical mechanics and invariance
of the frequency spectrum}
{Commun. Math. Phys., {\bf 87}, 1982, p. 365-383.}

\ref{[It]}{H. ITO}{Convergence of Birkhoff normal forms for 
integrable systems}
{Comment. Math. Helvetici, {\bf 64}, 1989, p. 412-461.}

\ref{[Ko]}{A.N. KOLMOGOROV}{Th\'eorie g\'en\'erale des syst\`emes 
dynamiques et m\'ecanique classique}
{Proc. Int. Congress of Math. 1954, vol. {\bf 1},
Amsterdam, 1957, p. 315-333.}

\ref{[Koz]}{V.V. KOZLOV}{Integrability and non-integrability in
Hamiltonian mechanics}
{Usp. Mat. Nauk. {\bf 38}, {\bf 1}, 1983, p.3-67. Russ. Math. Surv. 
{\bf 38 }, {\bf 1}, 1983, p. 1-76.}

\ref{[Li]}{A. LINDSTEDT}{Beitrag zur Integration der Differentialgleichungen
der St\"orungstheorie}
{Abh. K. Akad. Wiss. St. Petersburg, {\bf 31}, 4, 1882.}

\ref{[Ll]}{R. DE LA LLAVE}{On necessary and sufficient conditions
for uniform integrability of families of Hamiltonian systems}
{"International conference on dynamical systems, Montevideo 1995, 
Pitman Res. Notes. Math. Ser., {\bf 362}, Longman, Harlow, 1996, p.76-109.}

\ref{[Mo]}{J. MORALES RUIZ}{Differential Galois theory and non-integrability 
of Hamiltonian systems}
{Birkhauser, 1999.}

\ref{[Mo]}{J.K. MOSER}{On invariant curves of area 
preserving mappings of the annulus}
{Nachr. Akad. Wiss. G\"ottingen Math. Phys. Kl., {\bf 2 } , 1962,
p. 1-20.}

\ref{[Ne]}{N.N. Nehoroshev}{The behaviour of hamiltonian systems that 
are clase to integrable ones}
{Functional Analysis and its Applications, {\bf 5:4}, 1971; 
Usp. Mat. Nauk., {\bf 32:6}, 1977.}

\ref{[PM1]}{R. P\'EREZ MARCO}{Total convergence or general divergence
in Small Divisors} 
{Preprint, 2000.}

\ref{[PM2]}{R. P\'EREZ MARCO}{Linearization of holomorphic germs 
with resonant linear part} 
{Preprint, 2000.}

\ref{[Po]}{H. POINCAR\'E}{Les m\'ethodes nouvelles de la m\'ecanique
c\'eleste}{ Paris, 1892.}

\ref{[Ra]}{T. RANSFORD}{Potential theory in the complex plane}{London
Mathematical Society, Student Texts {\bf 28}, Cambridge University 
Press, 1995.}

\ref{[Rob]}{R.C. ROBINSON}{Generic properties of conservative
systems}{Amer. J. Math., {\bf 92}, 1970, p. 562-603, 
p. 897-906.}

\ref{[Ro]}{L.I. RONKIN}{Introduction to the theory of 
entire functions of several variables}{Translations of Mathematical
monographs, American Mathematical Society, {\bf 44}, 1974.}

\ref{[Ru1]}{H. R\"USSMANN}{Ueber das Verhalten analytischer Hamiltonscher 
Differentialgleichungen in der N\"ahe einer Glaichgewichtsl\"osung}
{Math. Ann. {\bf 154}, 1964, p. 285-300.}

\ref{[Ru2]}{H. R\"USSMANN}{\"Uber die Normalform analytischer 
Hamiltonscher Differentialgleichungen in der N\"ahe einer 
Gleichgewichtsl\"osung}
{Math. Ann. {\bf 169}, 1967, p. 55-72.}

\ref{[Ru3]}{H. R\"USSMANN}{On the convergence of power series 
transformations of analytic mappings near a fixed point}
{Preprint I.H.E.S., Paris, 1977.}

\ref{[Si1]}{C.L. SIEGEL}{On the integrals of canonical systems}{Annals
of Mathematics, {\bf 42}, {\bf 3}, 1941, p.806-822.}

\ref{[Si2]}{C.L. SIEGEL}{\"Uber die Existenz einer Normalform
analytischer Hamiltonscher Differentialgleichungen in der 
N\"ahe einer Gleichgewichtl\"osung}{Math. Annalen, {\bf 128},
1954, p.144-170.}

\ref{[Si-Mo]}{C.L. SIEGEL, J. MOSER}{Lectures on celestial mechanics}
{Springer-Verlag, {\bf 187}, 1971.}

\ref{[St1]}{L. STOLOVITCH}{Compl\`ete int\'egrabilit\'e singuli\`ere}
{C.R.A.S., {\bf 326}, 1998, p. 733-736.}

\ref{[St2]}{L. STOLOVITCH}{Normalisation holomorphe d'alg\`ebres 
de type Cartan de champs de vecteurs holomorphes singuliers}
{Pr\'epublications Universit\'e Paul Sabatier, {\bf 186}, mars 2000.}

\ref{[Wh]}{E.T. WHITTAKER}{On the solutions of dynamical problems in terms
of trigonometric series}{Proc. London Math. Soc., II 
{\bf 34}, 1902, p. 206-221.}

\ref{[Wi]}{A. WINTNER}{The analytic foundations of Celestial Mechanics}{
Princeton University Press, Princeton, New Jersey, 1941.}